\documentclass[a4paper,11pt]{article}
\usepackage{amsmath}
\usepackage{amsfonts}
\usepackage{amssymb}

\usepackage{latexsym}
\usepackage{graphicx}
\usepackage{extarrows}
\usepackage{epsfig}
\usepackage{mathrsfs}
\usepackage{latexsym}
\usepackage{graphicx}
\usepackage{epstopdf}
\usepackage{float}
\usepackage{caption}
\usepackage{setspace}
\usepackage{multirow}
\usepackage{algorithmicx}
\usepackage{algorithm}
\usepackage{algpseudocode}
\usepackage{authblk}

\newtheorem{theorem}{\bf Theorem}[section]
\newtheorem{lemma}[theorem]{\bf Lemma}

\newtheorem{remark}[theorem]{\bf Remark}

\newtheorem{defn}[theorem]{\bf Definition}

\def\D{\bf D}
\def\E{\bf E}

\def\eop{\hfill\rule{2.0mm}{2.0mm}}

\def\H{\mathbb{H}^2}

\def\C{\bf{C}}

\def\D{\bf D}
\def\E{\bf E}

\begin{document}

\title{\bf A sufficient condition for $n$-Best Kernel Approximation in Reproducing Kernel Hilbert Spaces}
\author[1]{Wei Qu}
\author[2]{Tao Qian \thanks{Corresponding author: tqian@must.edu.mo\\
 Funded by The Science and Technology Development Fund, Macau SAR (File no. 0123/2018/A3)}}
\author[3]{Guan-Tie Deng}
\affil[1]{School of Mathematical Sciences, Beijing Normal University, Beijing, China}
\affil[2]{Macau Center for Mathematical Sciences, Macau University of Science and Technology, Macau, China}
\affil[3]{School of Mathematical Sciences, Beijing Normal University, Beijing, China}


\maketitle
\begin{abstract}
 We show that if a reproducing kernel Hilbert space $H_K,$ consisting of functions defined on $\E,$ enjoys Double Boundary Vanishing Condition (DBVC) and Linear Independent Condition (LIC), then for any preset natural number $n,$ and any function $f\in H_K,$ there exists a set of $n$ parameterized multiple kernels ${\tilde{K}}_{w_1},\cdots,{\tilde{K}}_{w_n}, w_k\in {\bf E}, k=1,\cdots,n,$ and real (or complex) constants $c_1,\cdots,c_n,$ giving rise to a solution of the optimization problem
  \[ \|f-\sum_{k=1}^n c_k{\tilde{K}}_{w_k}\|=\inf \{\|f-\sum_{k=1}^n d_k{\tilde{K}}_{v_k}\|\ |\ v_k\in {\bf E}, d_k\in {\bf R}\ ({\rm or}\ {\bf C}), k=1,\cdots,n\}.\]
  By applying the theorem of this paper we show that the Hardy space and the Bergman space, as well as all the weighted Bergman spaces in the unit disc all possess $n$-best approximations. In the Hardy space case this gives a new proof of a classical result. Based on the obtained results we further prove existence of $n$-best spherical Poisson kernel approximation to functions of finite energy on the real-spheres.

\end{abstract}

MSC: {41A20; 41A65; 46E22; 30H20}

\smallskip

{\bf keywords:} Reproducing Kernel Hilbert Space, Double Boundary Vanishing Condition, $n$-Linearly Independent Condition, Hardy Space, Bergman space, Approximation by Rational functions of Certain Degrees
%

\bigskip


\section{Introduction}

Let $H$ be a complex Hilbert space consisting of functions defined in a topological space $\E$. Assume that the point evaluation functional $f(z)$ for any fixed $z\in \E$ is a bounded linear functional, i.e.,\[ |f(z)|\leq C_z \|f\|,\]
where $C_z$ is a constant depending on $z$. Then, according to the Riesz representation theorem there is a function $K_z$ with $z$ being a parameter such that
\[ f(z)=\langle f,K_z\rangle,\]
for all $z\in \E.$ In such case we say that $H$ is a reproducing kernel Hilbert space, abbreviated as RKHS, call $K_z$ the reproducing kernel of $H$. Denote by $H_K(\E)$ the Hilbert space $H$ whose corresponding reproducing kernel function is $K_z.$ Indeed, any RKHS can have only one reproducing kernel.
 A wide class of Hilbert spaces, including the classical Hardy ${\H}$-spaces, Bergman spaces, weighted Bergman spaces, and Sobolev spaces, etc., belong to the category of reproducing kernel Hilbert spaces (RKHSs). The subject $n$-best approximation in reproducing kernel Hilbert spaces include, as a particular case, the one called best approximations to Hardy space functions by rational functions of order not exceeding $n$. The present study amounts to extending the question and solving it in a wide class of Hilbert spaces.
  \par
In below we first provide an account of the related concepts in the classical Hardy space of the unit disc.
  Denote by ${\bf C}$ the complex plane and ${\bf D}$ the open unit disc in ${\bf C}.$ The Hardy space in the unit disc is defined, among other equivalent definitions,
 \[ \H({\bf D})=\{f:{\bf D}\to {\bf C}\ |\ f(z)=\sum_{k=0}^\infty c_kz^k, \sum_{k=0}^\infty |c_k|^2<\infty\}.\] It is a basic property of the Hardy space that for any $f\in \H({\bf D})$ there exists a boundary limit function, denoted $f(e^{it})\in L^2(\partial {\bf D}),$ in both the pointwise non-tangential limit sense and in the $L^2$-convergence sense as well.
 It is standard knowledge that under the inner product
 \[ \langle f,g\rangle =\frac{1}{2\pi}\int_0^{2\pi}f(e^{it})\overline{g}(e^{it})dt\]
 the space $\H(\D)$ forms a Hilbert space.

Of particular importance in the Hardy space theory are the functions
 $$k_w(e^{it})=\frac{1}{1-\overline{w}e^{it}}\quad {\rm and}\quad e_w^{\H}(z)=\frac{k_w}{\|k_w\|}=\frac{\sqrt{1-|w|^2}}{1-\overline{w}z}.$$
The function $k_w(z)$, where $w$ is considered as a parameter, is the reproducing kernel of the Hardy space $\H({\bf D}):$ By invoking the Cauchy formula it follows that for any $f\in \H({\bf D})$ there holds
\[ \langle f,k_w\rangle = f(w),\] and, subsequently,
\[ \langle f,\left(\frac{\partial}{\partial {\overline{w}}}\right)^lk_w\rangle = f^{(l)}(w), \qquad l=1,2,\cdots.\]

\begin{defn}
For any $n$ complex numbers $(w_1,\cdots,w_n)\in \D^n$ and $n$ complex numbers $(c_1,...,c_n)\in{\C}^n ,$  the function
\[ \sum_{k=1}^nc_kB_k(z)\] is called an $n$-Blaschke form, and an $n$-degenerate Blaschke form if $c_n\ne 0,$ where $\{B_k\}_{k=1}^n$ is the $n$-Takenaka-Malmquist (n-TM) system generated by the sequence $(w_1,\cdots,w_n),$
 \[ B_k(z)=\frac{\sqrt{1-|w_k|^2}}{1-\overline{w}_kz}\prod_{l=1}^{k-1}
 \frac{z-w_l}{1-\overline{w}_lz}.\]
\end{defn}
We note that the $n$-TM system $\{B_k\}_{k=1}^n$ is the orthonormalization of the $n$-system
\[ (\tilde{k}_{w_1},\cdots,\tilde{k}_{w_n}),\]
where
\begin{eqnarray}\label{df1} \tilde{k}_{w_k}(z)\triangleq\left(\frac{d}{d{\overline{w}}}\right)^{(l(w_k)-1)}
k_{w}(z)|_{w=w_k},\end{eqnarray}
called the \emph{multiple reproducing kernels}, where
\[l(w_k)\triangleq {\rm multiple \ number \ of }\ w_k \ {\rm in }\ (w_1,\cdots,w_k)\]
(\cite{Qian Wegert 2013,qian2013remarks}).
Besides the multiple reproducing kernels we also use \emph{normalized multiple
reproducing kernels}
\begin{eqnarray}\label{df2} \tilde{e}_{w_k}^{\H}(z)=\frac{\tilde{k}_{w_k}(z)}{\|\tilde{k}_{w_k}(z)\|}.\end{eqnarray}
 For fast expanding a given function into a TM system the \emph{adaptive Fourier decomposition}
(AFD) was proposed that is related to the Beurling-Lax decomposition of the Hardy space into the direct sum of the forward- and the backward-shift invariant subspaces (\cite{QWa,TQC}). AFD theory and algorithm have been generalized to matrix-valued functions defined in the disc (\cite{ACQS1}) and in the ball
 of several complex variables (\cite{ACQS2}).

  The $n$-best rational approximation problem in the Hardy space is formulated as follows. A pair of polynomials $(p,q)$ is said to be $n$-\emph{admissible} if it satisfies the following conditions: (i) $p$ and $q$ are co-prime;(ii) $q$ does not have zeros in ${\bf D};$ and (iii) the both degrees of $p$ and $q$ are at most $n$ (\cite{Qian Wegert 2013,Qcy}).\\

\noindent{\bf The $n$-best Rational Approximation Problem}: For $f\in \H({\bf D}),$ find an $n$-admissible pair of polynomials $(p_1,q_1)$ such that
\begin{eqnarray}\label{nbestrational} \|f-\frac{p_1}{q_1}\|_{\H({\bf D})}=\inf \{\|f-\frac{p}{q}\|_{\H({\bf D})} \ |\ (p,q) \ {\rm is\ }n{\rm - admissible}\}.\end{eqnarray}

The above optimization problem may be re-formulated as finding a non-degenerative Blaschke form
\[ \sum_{k=1}^n \langle f,B^{\bf w}_k\rangle B^{\bf w}_k(z),\] where the $B^{\bf w}_k$'s correspond to ${\bf w}=(w_1,\cdots,w_n),$
such that
\begin{eqnarray}\label{bestnBlaschke} \|f-\sum_{k=1}^n \langle f,B^{\bf w}_k\rangle B^{\bf w}_k\|=\inf \{\|f-\sum_{k=1}^n \langle f,B^{\bf v}_k\rangle B^{\bf v}_k\|_{\H({\bf D})} \ |\ {\bf v}=(v_1,\cdots,v_n) \in \D^n\}.\end{eqnarray}

Or, alternatively, we can ask the following question: Denotes by $f/{\rm span}\{\tilde{e}_{v_1}^{\H},\cdots,\tilde{e}_{v_{n}}^{\H}\}$ the projection of $f$ into the span of $\tilde{e}_{v_1}^{\H},\cdots,\tilde{e}_{v_{n}}^{\H}, (v_1,\cdots,v_n)\in \D^n.$ Find $(w_1,\cdots,w_n)\in \D^n$ such that
\[ \| f/{\rm span}\{\tilde{e}_{w_1}^{\H},\cdots,\tilde{e}_{w_{n}}^{\H}\}\| \]
is maximized over all $(v_1,\cdots,v_n)\in \D^n.$

There have been several proofs in the literature for existence of the above specified $n$-best rational approximation problem in the  classical Hardy spaces, see \cite{walsh1935interpolation} (J. L. Walsh, 1962), \cite{Ruck1978} (G. Buckebusch, 1978), \cite{Bara1986} (L. Baratchart), \cite{Qian Wegert 2013}, \cite{MQW2012}. In the last two articles the problem is reformulated in terms of $n$-Blaschke form.
In the Hardy space case practical algorithms, including the INRIA method (\cite{Baratchart1991}), cyclic AFD (\cite{Qcy}), and lately the gradient descent method in \cite{QWM}, can only claim to converge, in fact, to a local minimum.  A mathematical algorithm to find the global minimum, is now still being sought.

The present paper works with the RKHS context. In a general RKHS one has a set of analogous objects and can  raise the same $n$-best approximation question. Let $H_K$ be a reproducing kernel Hilbert space (RKHS) consisting of a class of functions defined in a topological space $\E,$ an open and connected set if it is in a larger topological space, with the reproducing kernel $K_w:$ $w\in {\bf E},$ that is, for any $f\in H_K,$
\[ \langle f,K_w\rangle=f(w).\]
We will also use the objects $\tilde{K}_w, \tilde{E}_w$ as, respectively, multiple reproducing kernel and multiple normalized reproducing kernel, similarly defined as
$\tilde{k}_a$ and $\tilde{e}^{\H}_a$ in, respectively, (\ref{df1}) and (\ref{df2}).

 For a fixed positive integer $n,$ the $n$-best question is formulated as follows: Find $n$ parameters ${\bf a}=(a_1,\cdots,a_n)\in \E^n$ that make the objective function
\begin{eqnarray}\label{simul} A(f;{\bf a})=\|f-\sum_{k=1}^n \langle f, B^{\bf a}_k\rangle_{H_K} B^{\bf a}_k \|_{H_K}\end{eqnarray}
minimized, where
\begin{eqnarray}\label{I want}
\sum_{k=1}^n \langle f, B^{\bf a}_k\rangle_{H_K} B^{\bf a}_k
\end{eqnarray}
 is called the $n$-\emph{kernel orthonormal form} of $f$ corresponding to the $n$-tuple
$(a_1,\cdots,a_n),$ where
 $(B^{\bf a}_1,\cdots,B^{\bf a}_k)$ is the G-S
 orthonormalization of the multiple kernels
 $(\tilde{E}_{a_1},\cdots,\tilde{E}_{a_k}), k\le n.$ Note that the above formulation is equivalent with the following minimization problem: Find $(a_1,\cdots,a_n)\in {\E}^n, (c_1,\cdots,c_n)\in {\C}^n,$ such that
 \begin{eqnarray}\label{problem}\|f-\sum_{k=1}^n c_k\tilde{K}_{a_k}\|=\inf \{\|f-\sum_{k=1}^n d_k\tilde{K}_{b_k}\|\ |\ b_k\in {\bf E}, d_k\in {\bf C}, k=1,\cdots,n\}.\end{eqnarray}
If for some $d_k$' and $b_k$'s $f=\sum_{k=1}^m d_k\tilde{K}_{b_k},$ then $f$ is said to be an \emph{$m$-kernel expansion}.

 We note that in the cases where the RKHS under study is the Hardy space inside the unit disc or the Hardy space in the upper-half complex plane, if $(B_1,\cdots,B_n)$ is the G-S orthonormalization of the $n$-tuple of the multiple reproducing kernels $(\tilde{e}_{a_1}^{\H},\dots,\tilde{e}_{a_n}^{\H}),$ then, by adding one more multiple reproducing kernel $\tilde{e}_{a_{n+1}}^{\H}$ to the $n$-sequence, the corresponding $(n+1)$-orthonormalization system, $(B_1,\cdots,B_n,B_{n+1}),$ is with the $(n+1)$-th term of the form
$B_{n+1}=\phi_n \tilde{e}_{a_{n+1}}^{\H},$ where $\phi_n$ is the Blaschke product, unique up to a uni-modular constant, defined by the first $n$ parameters $a_1,\cdots,a_n$ as its zeros, including the multiples. Indeed, TM systems are constructed in such way. In AFD, through a generalized backward shift operation, the TM systems are automatically generated (\cite{QWa}). It is a question whether there exist other types RKHSs that possess such or similar constructive property. From our observation it seems that only the Hardy spaces of the classical domains possess such property (see \cite{ACQS1,ACQS2}). In the weighted Bergman spaces of the classical domains this property does not hold (\cite{qu2018,qu2019}).

As technical preparation we need to recall the so called $\rho$-weak pre-orthogonal adaptive Fourier decomposition ($\rho$-Weak-POAFD) developed in the general RKHS context. Assume that ${\cal H}$ is a general Hilbert space with a dictionary parameterized by elements in $\E$, denoted $E_a, a\in \E$. Let $\rho\in (0,1).$  Suppose that we have obtained an $n$-term orthogonal expansion
\[ f=\sum_{k=1}^n \langle f,B_k\rangle B_k +g_{n+1},\]
where $(B_1,\cdots,B_n)$ is the G-S orthonormalization of a selected $n$-sequence $(a_1,\cdots,a_n),$ where the $a_k$'s are mutually different. Select $a_{n+1} ,$ different from all the already selected $a_k$'s, $k=1,\cdots, n,$ such that
\begin{eqnarray}\label{maximal} |\langle f,B^{a_{n+1}}_{n+1}\rangle|\ge \rho
\sup \{ |\langle f,B^{b}_{n+1}\rangle|\ |\ b\in \E\},\end{eqnarray}
where for any $b\in {\E}, (B_1,\cdots,B_n,B^b_{n+1})$ is the G-S orthonormalization of $(E_{a_1},\cdots,E_{a_n},E_b).$ Make such selections from $n=1$ and for all consecutive $n>1,$ we obtain $\rho$-Weak-POAFD (\cite{Q2D,qu2018}.

\begin{remark}  $\rho$-Weak-POAFD is available for all RKHSs. The selection criterion (\ref{maximal}) shows that it is a more optimal selection principle than the other types weak greedy algorithms in the classical literature (\cite{MaZ,LT}). When a dictionary satisfies BVC (see below), the selection corresponding to $\rho=1$ is available, called POAFD. POAFD has the optimal maximal selection at each algorithm step (\cite{Q2D}).  It reduces to AFD in the classical Hardy space (\cite{QWa}).
\end{remark}

In this paper we introduce what we call by \emph{Double Boundary Vanishing Condition} (DBVC) that will play an important role in the $n$-best optimization problem. Assume that the parameters set $\E$ is equipped with a topology. We used to work with the cases in which $\E$ is a region (open and connected) of the complex plane $\C$ or a region of the space of several complex variables ${\C}^n$ under its natural topology.  We now have the convention that, together with the finite boundary points, we add the infinite point, being included in the set of the boundary points if $\E$ is unbounded, that corresponds to the one-point-compactification of the original topological space. The added point is denoted $\infty.$ Taking ${\E}={\C}^+=\{{z}\in {\C}\ |\ {\rm Im}(z)>0\}$ as an example. $\E$ is equipped with the topology of $\C.$ By adding the $\infty$ point, the sequence of open sets $\{ z\in {\C}^+ \ |\ {\rm Im}(z)<\frac{1}{m}\ {\rm or} \ |z|>n\},$ were $m,n$ are positive integers, forms a basis of open neighborhoods of $\partial \E.$  A RKHS is said to satisfy DBVC if for any sequence $z_n\to \tilde{z}\in \overline{\E}$ and $w_n\to \tilde{w}\in \partial \E,$ and $z_n\neq w_n,$ there holds
\begin{eqnarray}\label{DBVC}
 \lim_{n\to \infty}\langle E_{z_n},E_{w_n}\rangle =0.\end{eqnarray}
 If DBVC holds, then we can show BVC (boundary vanishing condition) holds. That is, for any $f\in H_K$ and $ w_n\to \tilde{w}\in \partial \E,$ there holds
  \[ \lim_{n\to \infty} \langle f,E_{w_n}\rangle =0.\]
  We have the following
  \begin{lemma} \label{lemma} If $H_K$ is a RKHS satisfying DBVC, then it satisfies BVC.
  \end{lemma}
   \noindent{\bf Proof.} Let $f\in H_K.$ Since $H_K$ is a RKHS, by any type of the matching pursuit  algorithm, including POAFD and Weak-POAFD, one can find $(a_1,\cdots,a_n,\cdots),$ consisting of mutually different terms in ${\E},$ such that
  \[ f=\sum_{k=1}^\infty \langle f,B_k \rangle B_k,\]
  where for any $n,$ $(B_1,\cdots,B_n)$ is the G-S orthonormalization of some selected $(E_{a_1},\cdots,E_{a_n}), n=1,2,\cdots$ Then, for any $\epsilon>0,$ one can find a natural number $N$ such that
  \[\|f-\sum_{k=1}^N \langle f,B_k \rangle B_k\|\leq \epsilon/2.\]
  By invoking the Cauchy-Schwarz inequality we have
  \begin{eqnarray*}
  |\langle f,E_{w_n}\rangle|&=&|\langle f-\sum_{k=1}^N \langle f,B_k \rangle B_k,E_{w_n} \rangle| +|\langle \sum_{k=1}^N \langle f,B_k \rangle B_k,E_{w_n}\rangle| \\
  &\leq& \|f-\sum_{k=1}^N \langle f,B_k \rangle B_k\|+|\langle \sum_{k=1}^N \langle f,B_k \rangle B_k,E_{w_n}\rangle|\\
  &\leq& \epsilon/2 +|\langle \sum_{k=1}^N \langle f,B_k \rangle B_k,E_{w_n}\rangle|.
  \end{eqnarray*}
    We note that in the last summation the functions $B_1,\cdots,B_N$ can be expressed as linear combinations of $E_{a_1},\cdots,E_{a_N}.$ The inner products involving $B_1,\cdots,B_N$ then can be passed on to those with $E_{a_1},\cdots,E_{a_N},$ and thus DBVC can be used. As a result, the last term of the above inequality chain is less than $\epsilon$ if $n$ is large enough.  The proof is complete.\eop\\

  We need a condition on RKHS called \emph{n-Linearly Independent Condition} ($n$-LIC): If for a fixed $n$ and any mutually distinguish $w_1,\cdots,w_n$ the corresponding function set $\{E_{w_1},\cdots,E_{w_n}\}$ is linearly independent, then the RKHS is said to satisfy $n$-Linearly Independent Condition. This condition is rather mild, for, if it is not true, then a parameterized reproducing kernel is a linear expansion of some others. The latter implies that there exist $w_1,\cdots,w_{n},$ such that for all functions $f$ in the space there holds
  $f(w_n)=c_1f(w_1)+\cdots +c_{n-1}f(w_{n-1}),$ where $c_k$'s are fixed complex constants. A consequence of $n$-LIC, that is also the form that we use in the proof of our main Theorem \ref{main}, is that if $a_1,\cdots,a_k,b$ are mutually distinguish points in $\E,$ then the projection of $E_b$ into the span of $E_{a_1},\cdots,E_{a_k},$ or
  $\|E_b-\sum_{k=1}^k \langle E_b,B_k\rangle B_k\|=\sqrt{1-\sum_{l=1}^k |\langle E_b,B_l\rangle|^2}$ is nonzero, where $(B_1,\cdots,B_k)$ is the G-S orthonormalization of $(E_{a_1},\cdots,E_{a_k}), k\leq n-1.$

 The main result of this paper is
 \begin{theorem}
 A RKHS $H_K$ has a solution for the $n$-best optimization problem (\ref{simul}) in the open set $\E^n$ if the RKHS satisfies {\rm DBVC} and $n$-{\rm LIC}.
\end{theorem}

The precise statement of the theorem will be given in next section. The main effort of the proof is to show that under the conditions {\rm DBVC} and $n$-{\rm LIC} a solution exists and must situate in the open set $\E^n$ (interior solution). In both the theory (sifting process) and applications (model reduction) a solution being inside the open set is crucial, as having been seen in the complex Hardy space rational approximation theory (see, for instance, the enclosed references by Walsh, Baratchart, Qian, and Qu et al.). The main mechanism for such interior solutions is DBVC. In general RKHSs, a solution of the $n$-best may also happen at the boundary. Hence DBVC is not a necessary condition of existence of a general solution.

After proving the main theorem we verify that the weighted Bergman spaces in the disc satisfy DBVC and $n$-LIC, and thus conclude that the weighted Bergman spaces have $n$-best kernel approximations in the corresponding Hilbert space norms. Based on the obtained results we further prove existence of $n$-best spherical Poisson kernel approximation to functions of finite energy on the real-spheres. Except the classical Hardy spaces case, the other $n$-best existence results proved in this paper, including the version on RKHSs with a DBVC and $n$-LIC dictionary and the concrete examples with complex holomorphic function spaces and the spaces of functions of finite energy on the real-spheres, are all new results and proved for the first time.

\section{Existence of $n$-Best Approximation for RKHS with DBVC and $n$-LIC}

 \begin{theorem}\label{main} Let $H_K$ be a RKHS that satisfies DBVC  and $n$-LIC. Let $n$ be any but fixed positive integer. Then for any $f\in H_K,$ if $f$ by itself is not an $m$-kernel expansion form for $0\leq m<n,$ then there exists an $n$-tuple of parameters $(a_1,\cdots,a_n)\in \E^n,$ with $(B_1,\cdots,B_n)$ being the associated orthonormal systems such that
 \[ A(f;{\bf a})=\| f-\sum_{t=1}^n \langle f,B_t\rangle B_t\| \]
 attains the minimum value over all possible values arising from all the $n$-tuples in place of $(a_1,\cdots,a_n)$ in $\E^n.$
 \end{theorem}

 \noindent{\bf Proof of Theorem \ref{main}}. Denote ${\bf b}=(b_1,\cdots,b_n).$ It is obvious that $A(f;{\bf b})$ has a non-negative global infimum value for all ${\bf b}$ in ${\bf E}^n,$ call it $d.$ We show that this global infimum value is attainable at an interior point of ${\bf E}^n.$ Let ${\bf a}^{(k)}=(a^{(k)}_1,\cdots,a^{(k)}_n)$ be an $n$-tuple at which
 $A(f;{\bf a}^{(k)})<d+1/k.$ There then exists a subsequence tending to an $n$-tuple ${\bf a}$ in $\overline{\E}^n.$ Without loss of generality we can assume that the sequence ${\bf a}^{(k)}$ itself tends to ${\bf a}.$ We are to show ${\bf a}\in {\E}^n.$ Assume the opposite, which means that some coordinates of ${\bf a}$ are on $\partial{\E},$ and we will, in such case, introduce a contradiction. We divide the $n$ coordinates into two groups, $\mathbb{I}$ and $\mathbb{B},$ where for $l\in \mathbb{I}$ there holds
   $\lim_{k\to \infty}a^{(k)}_l=a_l\in \E;$ and for $l\in \mathbb{B}$ there holds
   $\lim_{k\to \infty}a^{(k)}_l=a_l\in \partial{\E}.$ We are assuming  $\mathbb{B}\ne \emptyset.$ Since $A(f;{\bf a}^{(k)})$ is the energy of $f$ onto the orthogonal complement of the span of the multiple reproducing kernels in the $n$-tuple ${\bf a}^{(k)},$ the energy quantity being irrelevant with the order of the elements in ${\bf a}^{(k)},$ we can assume, without loss of the generality, that the coordinates in $\mathbb{I}$ are all in front of those in $\mathbb{B}.$ The point is to show that, because $\lim_{k\to \infty}a^{(k)}_l\in \partial{\E}$ for $l\in \mathbb{B},$ the components $a^{(k)}_l$ of ${\bf a}^{(k)},$ if $l\in \mathbb{B},$ will have no contributions to the optimization of $A(f;{\bf a}).$
   To simplify the argument we may assume without loss of generality that for every $k$ the $n$-tuple ${\bf a}^{(k)}$ does not have multiple components, although the limiting $n$-tuple ${\bf a}$ may have. Let $l_0$ be the largest index for the indices in $\mathbb{I},$ then the indices $l_0+t, 0<t\leq n-l_0$ will be in the index range $\mathbb{B}.$ Since $\mathbb{B}\ne \emptyset,$ we have $l_0<n.$

   Let $R^{(k)}={\rm span}\{{E}_{a^{(k)}_{1}}, \cdots, {E}_{a^{(k)}_{n}}\}$ and $P^{(k)}$ the orthogonal projection to $R^{(k)};$ and likewise,
   $R^{(k)}_{\mathbb{I}}=$
   ${\rm span}\{{E}_{a^{(k)}_{1}}, \cdots, {E}_{a^{(k)}_{l_0}}\}$  and
   $P^{(k)}_{\mathbb I}$ the orthogonal projection mapping into $R^{(k)}_{\mathbb{I}}.$  It is easy to show that for the given function $f$, the projections $P^{(k)}_{\mathbb I}f$ have a
   limit as $k\to \infty,$ denoted $P_{\mathbb I}f,$ as the projection of $f$ into ${\rm span}\{\tilde{E}_{a_1},...,\tilde{E}_{a_{l_0}}\}.$ Denote $g=f-P_{\mathbb I}f.$

  The general form of the elements in the Gram-Schmidt orthonomalization of the system
   $\{E_{a^{(k)}_t}\}_{t=1}^n$ is
   \begin{eqnarray}\label{GS} {B}^{(k)}_{t}=\frac{{E}_{a^{(k)}_{t}}-
   \sum_{j=1}^{t-1}\langle {E}_{a^{(k)}_{t}},B_j^{(k)}\rangle B_j^{(k)}}{\|{E}_{a^{(k)}_{t}}-
   \sum_{j=1}^{t-1}\langle {E}_{a^{(k)}_{t}},B_j^{(k)}\rangle B_j^{(k)}\|},\end{eqnarray}
   where $\{B_1^{(k)},\cdots,B_{j}^{(k)}\}$ is the Gram-Schmidt orthonormalization of $\{{E}_{a^{(k)}_{1}},\cdots,{E}_{a^{(k)}_{j}}\}, 1\leq j\leq n.$

   We show that for any function $h$ in the reproducing kernel Hilbert space there holds
   \begin{eqnarray}\label{vanishing} \lim_{k\to \infty}\langle h,B_t^{(k)}\rangle =0,\quad\quad l_0<t\leq n.\end{eqnarray}
   Temporarily accepting (\ref{vanishing}), and using it for $h=f,$ we have
   \[ d^2=\lim_{k\to \infty}  (\|f\|^2-\sum_{t=1}^{l_0} |\langle f,B_t^{(k)}\rangle|^2) =\|f-Pf\|^2=\|g\|^2.\]
We note that $g\ne 0,$  for otherwise $f$ is an $m$-kernel form with $m=l_0<n,$ contrary with the assumption. $g\ne 0$ then implies $d>0.$
   Let $g^{(k)}_j :=f-\sum_{t=1}^j \langle f,B^{(k)}_t\rangle B^{(k)}_t, l_0\leq j\leq n.$ We have $\lim_{k\to \infty}g^{(k)}_j:=g_{l_0}=g, l_0\leq j\leq n.$ Find $a\in \E$ such that $|\langle g, E_a\rangle|=\delta>0.$ Let the new parameter matrix be
$$\qquad\qquad\qquad b_t^{(k)}=a_t^{(k)}, 1\leq t<n;$$
$$b_n^{(k)}=a,$$
where only the last column is different from the old. Then in the new system, using $\tilde{B}^{(k)}_t$ in place of $B_t^{(k)},\ 1\leq t \leq n,$ and $\tilde{P}^{(k)}$ in place of ${P}^{(k)},$ we have
\begin{eqnarray}\label{*}
  \lim_{k\to \infty}\|\tilde{P}^{(k)}f\|^2 &=&  \lim_{k\to \infty}
(\sum_{t=1}^{n-1} |\langle f,\tilde{B}^{(k)}_t\rangle|^2 +|\langle f,\tilde{B}^{(k)}_n\rangle|^2) \nonumber\\
   &=& \lim_{k\to \infty}(\|\sum_{t=1}^{l_0} |\langle f,{B}^{(k)}_t\rangle{B}^{(k)}_t \|^2+|\langle f,\tilde{B}^{(k)}_n\rangle|^2) \nonumber\\
  &=& \| P_{\mathbb I}f\|^2 + \lim_{k\to \infty}|\langle f,\tilde{B}^{(k)}_n\rangle|^2 \nonumber \\
  &=&\| Pf\|^2 + \lim_{k\to \infty}|\langle f,\tilde{B}^{(k)}_n\rangle|^2
\end{eqnarray}
where $\langle f,\tilde{B}^{(k)}_n\rangle=\frac{\langle f,(I-P^k_{ n-1})E_a\rangle}{\rho_k},$ where, as a consequence of LIC, $\rho_k \in (0,1].$ We further have $\langle f,(I-P^k_{ n-1})E_a\rangle=\langle(I-P^k_{ n-1})f,E_a\rangle,$ also $\lim_{k\to \infty}P^k_{ n-1}f=P_{l_0}f.$ Taking into account $(I-P_{l_0})f = g,$ and $\lim_{k\to \infty}\rho_k=\rho =\sqrt{1-\sum_{t=1}^{l_0}|\langle E_a,B_t\rangle|^2}\in (0,1],$ as a consequence of LIC again. The last equality chain (\ref{*}) finally equals
\begin{eqnarray*}
 \| \tilde{P}f\|^2=\| Pf\|^2+ |\frac{\langle g,E_a\rangle}{\rho}|^2 = \| f\|^2-\|g\|^2+ (\delta/\rho)^2
   > \|f\|^2-d^2.
\end{eqnarray*}
Or,
\[ \lim_{k\to \infty}\|f-\tilde{P}^{(k)}f\|^2< d^2,\]
being contrary with $d$ being the global infimum value of $A(f;{\bf b}), {\bf b}\in \E^n.$  The proof of the theorem is complete.

   Now we proceed to prove the relation (\ref{vanishing}) for $t=l_0+1,\cdots,n.$ First let $t=l_0+1.$ We have
   \begin{eqnarray*}
   \langle f,B^{(k)}_{l_0+1}\rangle = \langle f,\frac{E_{a^{(k)}_{l_0+1}}-\sum_{j=1}^{l_0}\langle E_{a^{(k)}_{l_0+1}},B^{(k)}_j\rangle B^{(k)}_j}{\sqrt{1-
   \sum_{j=1}^{l_0}|\langle E_{a^{(k)}_{l_0+1}},B^{(k)}_j\rangle|^2}}\rangle.
   \end{eqnarray*}

   Since $H_K$ satisfies DBVC,  from Lemma \ref{lemma},  $H_K$ also satisfies BVC. As a consequence,
   \begin{eqnarray}\label{first} \lim_{k\to \infty}\langle f,E_{a^{(k)}_{l_0+1}}\rangle =0.\end{eqnarray}

  Since $\lim_{k\to \infty}(a_1^{(k)},\cdots,a_{l_0}^{(k)})=(a_1,\cdots,a_{l_0})\in {\E}^{l_0},$ there exist the limits $\lim_{k\to \infty}B^{(k)}_j=B_j,$ being functions in $H_K,$ for $j=1,\cdots,l_0.$  Then BVC and the Cauchy-Schwarz inequality imply
   \begin{eqnarray}\label{second}|\langle E_{a^{(k)}_{l_0+1}},B^{(k)}_j\rangle|&=&|\langle E_{a^{(k)}_{l_0+1}},B_j\rangle| + |\langle E_{a^{(k)}_{l_0+1}},B^{(k)}_j-B_j\rangle|\nonumber \\
   &\leq & |\langle E_{a^{(k)}_{l_0+1}},B_j\rangle| + \|B^{(k)}_j-B_j\|\nonumber \\
   &\to&  0,\quad {\rm as}\ k\to \infty, \qquad j\le l_0.\end{eqnarray}
   In accordance with the relations (\ref{first}) and (\ref{second}), we have (\ref{vanishing}) for $t=l_0+1.$

   Now we prove (\ref{vanishing}) for $t>l_0+1.$ The induction hypotheses include that each term $B^{(k)}_j, 1\leq j\leq t-1,$ is a linear combination of $E_{a^{(k)}_s}, 1\leq s\leq j,$ while the coefficients of the linear combination are all constituted by sums and products between $\langle E_{a^{(k)}_s}, E_{a^{(k)}_{s'}}\rangle, 1\leq s', s\leq j,$ and divisions by $\sqrt{1-
   \sum_{l=1}^{s-1}|\langle E_{a^{(k)}_{s}},B^{(k)}_l\rangle|^2},
   1\leq s \leq j,$ without involving universal constants; and that the $$\lim_{k\to \infty}\langle f,B^{(k)}_l\rangle=0, \quad l_0<l\leq t-1.$$

   Write, in accordance with (\ref{GS}),
   \begin{eqnarray*}
   \langle f,B_t^{(k)}\rangle = \frac{1}{\|{E}_{a^{(k)}_{t}}-
   \sum_{j=1}^{t-1}\langle {E}_{a^{(k)}_{t}},B_j^{(k)}\rangle B_j^{(k)}\|}
   \left(\langle f,{E}_{a^{(k)}_{t}}\rangle -
   \sum_{j=1}^{t-1}\langle {E}_{a^{(k)}_{t}},B_j^{(k)}\rangle \langle f,B_j^{(k)}\rangle\right).\end{eqnarray*}
   The assumed DBVC, its consequence BVC, and the induction hypotheses together, establish
   \[ \lim_{k\to \infty}\langle f,{E}_{a^{(k)}_{t}}\rangle =0, \
   \lim_{k\to \infty}\langle {E}_{a^{(k)}_{t}},B_j^{(k)}\rangle =0, \ {\rm and}
   \ \lim_{k\to \infty}\langle f,B_j^{(k)}\rangle=0, \ j\leq t-1.\]
   Therefore,
   \[ \lim_{k\to \infty}\langle f,B_t^{(k)}\rangle =0.\]
   Based on the mathematical induction principle the proof is complete.\eop

   \begin{remark}
   A large amount commonly used Hilbert spaces are RKHSs in which DBVC and LIC are satisfied.  The above theorem guarantees that such RKHSs have $n$-best kernel approximations. The recently developed cyclic and gradient descent algorithms (\cite{Qcy,QWM}) for Hardy spaces are adaptable to abstract RKHSs with DBVC and LIC. The proof of the existence result guarantees convergence of the adapted algorithms in abstract spaces. It, in particular, serves as a useful reference in learning theory for simultaneously selecting $n$-parameters to optimize an energy-based objective function.
   \end{remark}



\tabcolsep 30pt
\vspace*{7pt}

\end{document}